\documentclass{amsart}

\newtheorem{thm}{Theorem}
\newtheorem{prop}{Proposition}
\newtheorem{lemma}{Lemma}

\newcommand{\con}{\equiv}

\newcommand{\ndiv}{\nmid}
\newcommand{\modd}[1]{\; ( \text{mod} \; #1)}
\newcommand{\bstack}[2]{#1 \atop #2}
\newcommand{\tstack}[3]{#1\atop {#2 \atop #3}}

\newcommand{\maps}{\rightarrow}

\newcommand{\y}{\mathbf{y}}

\newcommand{\al}{\alpha}

\newcommand{\ga}{\gamma}
\newcommand{\del}{\delta}
\newcommand{\ep}{\epsilon}

\newcommand{\om}{\omega}

\newcommand{\Sig}{\Sigma}

\newcommand{\A}{\mathcal{A}}
\newcommand{\Bcal}{\mathcal{B}}

\newcommand{\Hcal}{\mathcal{H}}

\newcommand{\U}{\mathcal{U}}
\newcommand{\Ucal}{\mathcal{U}}

\newcommand{\V}{\mathcal{V}}
\newcommand{\Vcal}{\mathcal{V}}

\newcommand{\h}{{\bf h}}
\newcommand{\hbf}{{\bf h}}

\newcommand{\mbf}{{\bf m}}

\newcommand{\tbf}{{\bf t}}
\newcommand{\ubf}{{\bf u}}
\newcommand{\vbf}{{\bf v}}
\newcommand{\x}{{\bf x}}
\newcommand{\xbf}{{\bf x}}
\newcommand{\ybf}{{\bf y}}
\newcommand{\zbf}{{\bf z}}

\newcommand{\PP}{\mathbb{P}}

\newcommand{\C}{\mathbb{C}}
\newcommand{\F}{\mathbb{F}}

\newcommand{\Q}{\mathbb{Q}}
\newcommand{\R}{\mathbb{R}}

\newcommand{\Z}{\mathbb{Z}}

\newcommand{\beq}{\begin{equation}}
\newcommand{\eeq}{\end{equation}}
\usepackage{amsmath,amssymb,enumerate}
\usepackage{color}

\begin{document}
\title{Counting rational points on smooth cyclic covers}
\author{ D. R. Heath-Brown}
\address{Mathematical Institute, 24-29 St. Giles', Oxford, OX1 3LB}
\email{rhb@maths.ox.ac.uk}

\author{Lillian B. Pierce}
\address{Mathematical Institute, 24-29 St. Giles', Oxford, OX1 3LB}
\email{lillian.pierce@maths.ox.ac.uk}

\subjclass{14G05 (11D45, 11L40)}

\begin{abstract}
A conjecture of Serre concerns the number of rational points of bounded height on a finite cover of projective space $\PP^{n-1}$. In this paper, we achieve Serre's conjecture in the special case of smooth cyclic covers of any degree when $n \geq 10$, and surpass it for covers of degree $r \geq 3$ when $n > 10$. This is achieved by a new bound for the number of perfect $r$-th power values of a polynomial with nonsingular leading form, obtained via a combination of an $r$-th power sieve and the $q$-analogue of van der Corput's method.
\end{abstract}

\maketitle

\section{Introduction}

Let $F(\xbf) \in \Z[x_1,...,x_n]$ be an irreducible form
of degree $mr$ with $r \geq 2$, $m \geq 1$, such that the projective
hypersurface defined by $F({\bf x})=0$ is smooth. In this paper we
will investigate the number of integer solutions to \beq\label{power_eqn}
y^r = F(x_1, \ldots, x_n)
\eeq
with $|x_i| \leq B$. Our interest stems from the fact that an upper bound for the number of
such points provides an upper bound for the number of rational points
on cyclic covers of $\mathbb{P}^{n-1}$.

The density of rational points on covers of projective space is the
subject of a well-known conjecture of Serre. Precisely, given a finite cover $\phi: X \maps \mathbb{P}^{n-1}$ over
$\Q$, where $n\ge 2$, define the counting function \[ N_B(\phi) = \#\{P \in X(\Q): H(\phi(P)) \leq B \}.\]
Here $H$ is the usual multiplicative height function on $\mathbb{P}^{n-1}$.
Using a sieve method, Serre \cite{Ser89} proved that there exists $\ga<1$ such that
\beq\label{Ser_res}
N_B(\phi) \ll B^{(n-1)+\frac{1}{2}}(\log B)^\ga,
\eeq
as long as the degree of $f$ is at least two. In fact,
however, Serre conjectures that \beq\label{SConj}
N_B(\phi) \ll B^{n-1}(\log B)^c,
\eeq
for some $c$, for covers of any degree $r\ge 2$ (see Theorems 3, 4 of Chapter 13 in \cite{Ser89}).

Several results are known in this direction. Broberg \cite{Bro03} has
applied Heath-Brown's determinant method \cite{HB02} to prove results
for covers of $\PP^1$ and $\PP^2$. In the case of $\PP^1$, Broberg
proves that for $\phi: X \maps \PP^1$ of degree $r \geq 2$, \[ N_B(\phi) \ll_{\phi,\ep} B^{2/r + \ep}.\]
For $\phi: X \maps \PP^2$ of degree $r \geq 3$, Broberg proves that \[ N_B(\phi) \ll_{\phi,\ep} B^{2 + \ep},\]
and if $\phi$ is of degree $2$, then
\[ N_B(\phi) \ll_{\phi,\ep} B^{9/4+\ep}.\]
These results nearly prove Serre's conjecture (\ref{SConj}) for $n=2,3$.

Recently, Munshi \cite{Mun09} considered the case in which $\phi$ is a
\emph{smooth cyclic} cover of $\PP^{n-1}$, given by an equation of the type
(\ref{power_eqn}) with $F$ a nonsingular form.  In this situation he proves that for all $n \geq 2$, one has
\beq\label{MB}
N_B(\phi) \ll_\phi B^{n-\frac{n}{n+1}} (\log B)^{\frac{n}{n+1}}. \eeq
Note that if $n \geq 2$, this improves on (\ref{Ser_res}), and even
approaches Serre's conjecture in the limit as $n \maps \infty$.

Unpublished results of Salberger on the dimension growth conjecture \cite[Conjecture 2]{HB02} imply the truth of Serre's conjecture with $(\log B)^c$ replaced by $B^\ep$ for covers $\phi$ given by (\ref{power_eqn}) with $F$ a form of precisely degree $r$, for any $r \geq 2$. 
But note that while (\ref{SConj}) is as
good as one can hope for in complete generality, one might expect an
estimate of the shape
\[ N_B(\phi) \ll B^{n-m(r-1)}(\log B)^c\]
for covers of Munshi's type, under favourable circumstances. Thus
Serre's conjecture probably does not reflect the whole truth in this
area.

In this paper, we again start from the foundation of Munshi's
approach: given a form $F$ as above, the equation (\ref{power_eqn})
defines a variety $X$ in weighted projective space
$\mathbb{P}(m,1,\ldots, 1)$, where the first coordinate $y$ has weight
$m$ and the coordinates $x_1, \ldots, x_n$ have weight 1. The variety
$X$ can now be regarded as a cyclic $r$-sheeted cover of
$\mathbb{P}^{n-1}$, given explicitly by the map $\phi: X \maps
\mathbb{P}^{n-1}$ that takes the point $(y,x_1, \ldots, x_n)$ to
$(x_1, \ldots, x_n)$. Thus our attention turns to counting perfect $r$-th power values of
the form $F(x_1, \ldots, x_n)$ with $|x_i| \leq B$.

For convenience, we employ a smooth non-negative weight $w: \R^n \maps
\R$ such that $w \geq 1$ on the unit cube $[-1/2,1/2]^n$, is supported in
$[-1,1]^n$, and satisfies the differential inequality \[| \frac{\partial^\al}{\partial {\bf x}^\al} w({\bf x})| \ll_\al 1,\]
for all multi-indices $\al=(\al_1, \ldots, \al_n)$. Define the
normalized weight function \beq\label{wt_dfn}
w_B({\bf x}) = w({\bf x}/B),
\eeq
so that $w_B$ is supported in the box $\Bcal=[-B,B]^n$ and satisfies
the differential inequalities \[| \frac{\partial^\al}{\partial {\bf x}^\al} w({\bf x})| \ll_\al B^{-|\al|},\]
where $|\al| = \al_1 + \cdots + \al_n.$
We then define the counting function \beq\label{N_dfn}
N_{w,B} (F) = \sum_{y \in \Z} \sum_{\bstack{{\bf x} \in
   \Z^n}{F({\bf x}) = y^r}} w_B({\bf x}) ,
   \eeq
with the aim of proving upper bounds of the form \[ N_{w,B}(F) \ll B^{n-\del},\]
for some $\del>0$ independent of the degree of $F$.

The counting function $N_B(\phi)$ associated to the cyclic cover $\phi:X
\maps \PP^{n-1}$, where $X$ is defined by (\ref{power_eqn}), satisfies
the relation \[ N_B(\phi) \leq N_{w,2B}(F).\]
Munshi employed the square sieve, adapted to count perfect $r$-th powers, in
order to count solutions to (\ref{power_eqn}). This method ultimately
requires one to estimate mixed character sums, for which bounds Munshi
employed bounds of Deligne \cite{Del74} and Katz \cite{Kat02},
\cite{Kat07}. We also use a version of the power sieve, but we sieve
over certain almost-primes, instead of primes, and this allows us to apply
the $q$-analogue of van der Corput's method; a similar combination has been
used previously in \cite{Pie06}, \cite{HB08}.

With no extra effort we can handle equations of the form
(\ref{power_eqn}) in which $F$ is a general polynomial $f$, not
necessarily homogeneous, whose leading form is nonsingular.
Ultimately, we prove the following theorem for the counting function $N_{w,B}(f)$ defined as in (\ref{N_dfn}) for solutions to the equation \beq\label{power_eqn2}
y^r = f(x_1, \ldots, x_n).
\eeq
\begin{thm}\label{form_thm}
Let $f(\xbf) \in \Z[x_1,\ldots, x_n]$ be a polynomial of degree $d\ge
3$, and assume that its leading form is nonsingular. Then for any $r \geq 2$, the counting function $N_{w,B}(f)$ for the number of solutions to (\ref{power_eqn2}) satisfies 
\[N_{w,B}(f) \ll \left\{\begin{array}{cc} B^{n-3n/(2n+10)}(\log B)^2, & n\ge 8,\\
B^{n-n(n-2)/(6n+4)}(\log B)^2, & 2\le n\le 8,\end{array}\right.\]
where the implied constant depends on $f,d,n.$ \end{thm}
The reader should recall that if $f$ is a polynomial of degree $d$
then its leading form is defined to be the form composed of all those
terms in $f$ with degree exactly $d$.

This theorem fails to hold for $d=2$. That the method of proof
fails to hold is visible from the inapplicability of Proposition
\ref{Lemma_Poisson_fg}, which requires $\deg f \geq 3$. But in fact
the statement of the theorem is also false for $d=2$ and $n > 10$,
since it is well known that there are $\gg B^{n-1}$ values $\x\in
[-B,B]^n$ for which $x_1^2+\ldots+x_n^2$ is a square. (This follows
from Theorems 5, 6 and 8 of Heath-Brown \cite{HB96}, for example.)

As an immediate corollary, we have:
\begin{thm}\label{proj_thm}
Let $\phi : X \maps \PP^{n-1}$ be a smooth finite cyclic cover given by the equation (\ref{power_eqn}) with $F$ a nonsingular form of degree $mr$. Suppose
either that $\phi$ has degree $r \ge 3$, or that $r=2$ and $m \geq 2$.  Then
\[ N_B(\phi) \ll \left\{\begin{array}{cc} B^{n-3n/(2n+10)}(\log B)^2, & n\ge 8,\\
B^{n-n(n-2)/(6n+4)}(\log B)^2, & 2\le n\le 8.\end{array}\right.\]
\end{thm}

In the case $r=2$, $m=1$, Theorem \ref{form_thm} no longer applies directly, but in this case the function $F$ in (\ref{power_eqn}) is a
quadratic form, and it is well known that $N_{w,B}(F)\ll_F
B^{n-1}\log B$.  (See Theorems 5--8 of Heath-Brown \cite{HB96}, for example.)
Assembling this with the relevant result of Theorem \ref{proj_thm} for $r=2$, we therefore have the following result for all smooth finite covers of degree $2$ given by (\ref{power_eqn}):
\begin{thm}\label{proj_thm2}
Let $\phi : X \maps \PP^{n-1}$ be a smooth finite cover of degree $2$,
given by the equation (\ref{power_eqn}) with $F$ nonsingular. Then \[ N_B(\phi)  \ll \left\{\begin{array}{cc} B^{n-1}(\log B)^2, & n\ge 10,\\
B^{9-27/28}(\log B)^2, & n=9,\\
B^{n-n(n-2)/(6n+4)}(\log B)^2, & 2\le n\le 8. \end{array}\right.\]
\end{thm}

We therefore see that, for cyclic covers of any degree with $F$ nonsingular, we can
achieve Serre's conjecture (\ref{SConj}) for $n\ge 10$,
and indeed surpass it for $n>10$ and degree $r\ge 3$. Moreover we improve on Munshi's
bound (\ref{MB}) for $n\ge 8$.

There is some prospect of a better result if one could treat the sum $T(\hbf)$ in (\ref{T_triple}) without splitting
into residue classes $k\modd{q_1}$.  If this were possible we would
expect  a gain of
$q_1^{1/2}$ at this stage.  What would be required is
an estimate of the shape
\[\sum_{\x\modd{p}}\chi_1(f(\x)+g(\x))\chi_2(f(\x))e_p(\mathbf{c}\cdot\x)
\ll p^{n/2}\]
in which $f(\x)$ and $g(\x)$ are polynomials in $n$ variables, having
smooth leading forms, and in which ${\rm deg}(f)>{\rm deg}(g)\ge 2$.
Katz \cite{Kat07} proves related results, but his theorems appear not to
cover the case required here.

\section{The power sieve}
We begin by formulating the sieve inequality we will use to count
integer solutions to (\ref{power_eqn2}). The method has its origins in
the ``square-sieve'' of Heath-Brown \cite{HB84}.
Following the presentation of
Munshi, we define a character that detects perfect $r$-th powers, in
analogy to the Legendre symbol used to detect perfect squares. For any prime $p$, since $\mathbb{F}_p^*$ is a cyclic group,
there is a non-canonical isomorphism \[ \theta_p : \mathbb{F}_p^* \maps \mu_{p-1} \]
onto the set $\mu_{p-1}$ of $(p-1)$-th roots of unity in $\C^*$. For
each $p \con 1 \modd{r}$, fix such a $\theta_p$. On the other hand, for such $p$, for every element $a \in \mathbb{F}_p^*$, the
quantity $a^{\frac{p-1}{r}}$ is a well-defined $r$-th root of
unity. Thus we can define a primitive Dirichlet character modulo $p$
by setting \[ \chi_p(n) = \theta_p(\bar{n}^{\frac{p-1}{r}})\]
for $(n,p)=1$ and $\chi_p(n)=0$ for $p|n$. Note that if $n$ is such
that $(n,p)=1$ and $n=m^r$ for some $m$, then \[ \chi_p(n) = \chi_p((\bar{m}^r)^{\frac{p-1}{r}}) = 1,\]
so that $\chi_p$ detects $r$-th powers, as desired (with the
possibility of over-counting). We will require characters to composite moduli, so for $q=p_1\cdots
p_k$ with primes $p_1, p_2, \ldots, p_k \con 1 \modd{r}$, we
define \[ \chi_{q}(n) = \chi_{p_1}(n) \chi_{p_2}(n) \cdots \chi_{p_k}(n).\]
Then it still remains true that for $(n,q)=1$ such that $n=m^r$ we
have $\chi_{q}(n)=1$, although now we may be over-counting $r$-th
powers even more significantly.

We now describe an $r$-th power sieve using the characters $\chi_q$; we
have specialized the statement of the following lemma to suit our
particular needs, but a more flexible formulation may be found in
\cite{Pie06} (stated there for the case $r=2$, but easily generalized
to $r \geq 2$). \begin{lemma}\label{sieve_lemma}
Let $\chi_q$ be the multiplicative character modulo $q$ defined as above.
Let $\A = \{uv:  u \in \U, v \in \V \}$ where $\U$ and $\V$ are disjoint sets of primes satisfying $p \con 1 \modd{r}$. Let $A =
\# \A$, $U = \# \U,$ and $V = \# \V$, so that $A=UV$. Furthermore,
assume that $V^3 \ll A$. Let $\om$ be a non-negative weight such that $\om(n)=0$ for $|n| \geq \exp( \min (U,V))$. Then

\begin{eqnarray} \sum_{n\not=0} \om (n^{r})
	&  \ll& A^{-1} \sum_n \om(n) + A^{-2} \sum_{v,v' \in \V}
       \sum_{u \neq u' \in \U}
		\left| \sum_n \om(n) \chi_{uv}(n)
                 \overline{\chi_{u'v'}(n)}\right|  \nonumber \\
	& & + \; UA^{-2} \sum_{v \neq v' \in \V} \left| \sum_{n}
         \om(n) \chi_v(n) \overline{\chi_{v'}(n)}
       \right|. \label{sig_lemma_eqn} \end{eqnarray}
\end{lemma}

We will refer to the terms on the right hand side of the $r$-th power
sieve (\ref{sig_lemma_eqn}) respectively as the trivial leading term, the
main sieve, and the prime sieve. To prove the lemma, consider
\[ \Sig = \sum_n \om(n) \left| \sum_{q \in \A} \chi_q(n) \right|^2.\]
Each $n$ is summed with non-negative weight, and in particular, if
$n=m^r\not=0$ and $\om(n) \neq 0$, then \[ \sum_{q \in \A} \chi_q(n) = \sum_{q \in \A} \chi_q(m^r) =
\sum_{\bstack{q \in \A}{ (q,m)=1}} 1
	\geq A - \sum_{\bstack{q \in \A}{(q,m) \neq 1}} 1 \gg A.\]
The last step follows since
$\om(n)$ is nonzero only if $|n| < \exp( \min(U, V))$, so that $\nu(m)
\ll \min(U,V)=o(A),$ where $\nu(m)$ denotes the number of distinct
prime divisors of $m$. Thus \beq\label{Sig_to_A}
\Sig \gg A^2\sum_{n\not=0} \om(n^r).
\eeq
But also
\begin{eqnarray} \Sig & = & \sum_{q,q' \in \A} \sum_n \om(n) \chi_{q}(n)
\overline{\chi_{q'}(n)} \nonumber\\
	& = & \sum_{q  \in \A} \sum_{n} \om(n) \chi_q(n)\overline{\chi_q(n)}
		+ \sum_{\bstack{q \neq q' \in \A}{(q,q')=1}} \sum_n
               \om(n) \chi_{q}(n) \overline{\chi_{q'}(n)}
               \label{sig_sums} \\
	& & +\; \sum_{\bstack{q\neq q' \in \A}{(q,q')\neq 1}} \sum_n
       \om(n) \chi_{q}(n) \overline{\chi_{q'}(n)} \nonumber. \end{eqnarray}
The first term in (\ref{sig_sums}) is bounded above by $A
\sum_n\om(n)$. The second term in (\ref{sig_sums}) will belong to the
main term in the sieve. The last term in (\ref{sig_sums}) may be broken into two subsums
$S(\U) + S(\V)$, where
\[S(\U)=   \sum_{v \in \V} \sum_{u \neq u' \in \U} \sum_{n} \om(n)
\chi_{uv}(n) \overline{\chi_{u'v}(n)}\] and
\[S(\V)= \sum_{u \in \U} \sum_{v \neq v' \in \V} \sum_{n} \om(n)
\chi_{uv}(n) \overline{\chi_{uv'}(n)}. \]
The sum $S(\U)$ is simply included in the main sieve term, but $S(\V)$ requires a different approach. We split it into two
further pieces, writing: \begin{eqnarray*}
S(\V) &=&  \sum_{u \in \U} \sum_{v \neq v' \in \V} \sum_{\bstack{n}{u
   \ndiv n}} \om(n) \chi_{v}(n) \overline{\chi_{v'}(n)}\\
	&=& U\sum_{v \neq v' \in \V} \sum_{n} \om(n) \chi_v(n)
       \overline{\chi_{v'}(n)}
	- \sum_{u \in \U} \sum_{v \neq v' \in \V} \sum_{\bstack{n}{u
           |n}} \om(n)\chi_v(n) \overline{\chi_{v'}(n)}\\
	&=& M(\V) - E(\V),
	\end{eqnarray*}
say. The term $M(\V)$ now gives rise to the third term in
(\ref{sig_lemma_eqn}), the prime sieve. The
error term $E(\V)$ can be bounded above in absolute value to give
\[ |E(\V)| \leq \sum_{u \in \U} \sum_{v \neq v' \in \V}
\sum_{\bstack{n\not=0}{u |n}} \om(n)
	 \leq V^2 \sum_{n\not=0} \om(n)\nu(n),\]
where as usual $\nu(n)$ denotes the number of distinct prime divisors
of $n$. By assumption, if $\om(n) \neq 0$ for some $n\not=0$ then, $\nu(n) \ll
\min(U,V)$. Thus \[ |E(\V)| \ll V^3\sum_n \om (n),\]
which is dominated by the trivial leading term as long as $V^3 \ll A$.
Thus, under this assumption, we have shown that \begin{eqnarray*}
|\Sig| & \ll & A \sum_{n} \om(n) + \sum_{v,v' \in \V} \sum_{u \neq u' \in \U}
		\left| \sum_n \om(n) \chi_{uv}(n)
                 \overline{\chi_{u'v'}(n)}\right|  \\
	& & \hspace{2cm}\mbox{}+  \; U\sum_{v \neq v' \in \V} \left|
         \sum_{n} \om(n)
         \chi_v(n) \overline{\chi_{v'}(n)} \right|. \end{eqnarray*}
The result of the lemma then follows by comparison with (\ref{Sig_to_A}).
\bigskip

We will apply the $r$-th power sieve using the sets \beq\label{Udef}
\Ucal = \{ \text{primes} \; u \con 1 \modd{r} :  Q^\al < u
\leq 2Q^\al \}, \eeq
\beq\label{Vdef}
\Vcal =\{ \text{primes} \; v \con 1 \modd{r} : Q^{1-\al} < v
\leq 2Q^{1-\al} \}, \eeq
where $Q=B^{\del}$ for some $\del>0$, and the exponent $\al$ is a real
parameter satisfying $2/3 \leq \al<1$; these parameters will be chosen later. Note in particular that
under these conditions, $V^3 \leq A$, and we may assume that the
sieving primes $u$ and $v$ do not divide $r$.

Recall the smooth weight function $w_B(\x) = w (\x/B)$ given in
(\ref{wt_dfn}). We will define the sieve weight by \[ \om(n) = \sum_{\bstack{\x \in \Z^n}{f(\x)=n}} w_B(\x).\]
Then \beq\label{N_upper_bd}
N_{w,B} (f) \ll \om(0) + \sum_{n\neq 0} \om (n^r).
\eeq
We need to handle separately the contribution to $N_{w,B}(f)$ arising from terms
with $f(\x)=0$.  There are many estimates in the literature covering
this situation.  For example Heath-Brown \cite[Theorem 2]{HB94} gives
a bound $O(B^{n-3+15/(n+5)})$, which is adequate for Theorem \ref{form_thm}.

For the remainder of (\ref{N_upper_bd}), the leading term in the sieve has upper bound
\beq\label{lead_term}
A^{-1}\sum_n \om(n) \ll Q^{-1}(\log Q)^2 \sum_{\x} w_B(\x)  \ll
B^nQ^{-1}(\log Q)^2 .
\eeq
If we take $Q=B^\del$, we therefore see that this contributes $O(B^{n-\del}(\log B)^2)$ to $N_{w,B}(f)$.

Our principal task is to estimate the main sieve term, namely \begin{eqnarray}\label{T_dfn}
\sum_n \om (n) \chi_{uv}(n) \overline{\chi_{u'v'}(n)}
	&=& \sum_{\x \in \Z^n} w_B(\x) \chi_{uv}(f(\x))
       \overline{\chi_{u'v'}(f(\x))}\nonumber \\ &= &\sum_{\x \in \Z^n} w_B(\x) \chi_{q_1}^*(f(\x)) \chi_{q_2}^*(f(\x)),
	\end{eqnarray}
where for convenience we have defined $q_1=uu'$ to be the product of
the ``large'' primes $u,u'$, and $q_2=vv'$ to be the product of the
``small'' primes.  Moreover we have set
\beq\label{chi_dfn}
\chi^*_{q_1}(n)=\chi_{u}(n) \overline{\chi_{u'}(n)}\qquad\mbox{and}
\qquad \chi^*_{q_2}(n)=\chi_{v}(n) \overline{\chi_{v'}(n)}. \eeq

In fact, for the prime sieve term we shall need to consider a similar
sum with $q_1$ being prime.  We therefore prove the following more general
result for weighted character sums of the form (\ref{T_dfn}). \begin{prop}\label{T_prop}
Let $q_1$ and $q_2$ be coprime integers and suppose that $q_1$ is either
prime or a product $p_1p_2$ of primes satisfying
$p_1<p_2<2p_1$. Write $p=q_1$ if $q_1$ is prime, or $p=p_1$
if $q_1=p_1p_2$.
Let $\chi_{q_1}$ and $\chi_{q_2}$ be multiplicative characters modulo
$q_1$ and $q_2$ respectively, and suppose that $\chi_{q_1}$ is non-principal.  Define \beq\label{Tq_dfn}
T(q_1,q_2)= \sum_{\x \in \Z^n} w_B(\x) \chi_{q_1}(f(\x)) \chi_{q_2}(f(\x)),
\eeq
where $B\ge q_2$.  Then \beq\label{T_prop_eqn}
T(q_1,q_2)\ll_{f} B^{n/2}q_1^{1/2}q_2^{n/2} + B^{n/2}q_1^{(n+2)/4}
+B^np^{-(n-2)/4}.
\eeq
\end{prop}
We will apply this result to the main sieve with $q_1,q_2$ each being
a product of two primes (not necessarily distinct in the case of
$q_2$), and to the prime sieve with $q_1,q_2$ being distinct primes. Theorem \ref{form_thm} will then follow, as we will show in \S \ref{bst}.

\section{The $q$-analogue of van der Corput's method}

We begin our proof of Proposition \ref{T_prop} by applying the $q$-analogue of van
der Corput's method.  The sum $T(q_1,q_2)$ involves a character to
modulus $q_1q_2$, and the effect of our version of van der Corput's
method is to produce a sum involving a character with a smaller
modulus, namely $q_1$. To do this we take $H=[B/q_2]$
and let $\Hcal$ denote the set of integer $n$-tuples in $[1,H]^n$,
so that $\#\Hcal = H^n$. Then \begin{eqnarray*}
H^n T(q_1,q_2) &=& \sum_{\h \in \Hcal} \sum_{\x} w_B(\x+ q_2\h) 
\chi_{q_1}(f(\x+ q_2\h)) \chi_{q_2}(f(\x+ q_2\h )) \\
	&=&  \sum_{\x \in [-B-Hq_2,B-q_2]^n}
       \chi_{q_2}(f(\x))  \sum_{\h \in \Hcal} w_B(\x+ q_2\h
       )\chi_{q_1}(f(\x+ q_2\h)).
	\end{eqnarray*}
Applying Cauchy-Schwarz, \beq\label{H2T2}
H^{2n} |T(q_1,q_2)|^2 \leq \Sig_1 \Sig_2,
\eeq
where \begin{eqnarray*}
\Sig_1 & = & \sum_{\x \in [-B-Hq_2,B-q_2]^n} | \chi_{q_2}(f(\x))|^2,\\
\Sig_2 & = &  \sum_{\bstack{\x}{(f(\x),q_2)=1}} | \sum_{\h \in \Hcal}
w_B(\x+ q_2\h)\chi_{q_1}(f(\x+ q_2\h ))|^2. \end{eqnarray*}
It will be convenient to drop the condition that $(f(\x),q_2)=1$ in
$\Sig_2$; by positivity this will still produce an upper bound. Then,
expanding the resulting sums in $\Sigma_2$, we have \[\Sig_2 \leq   \sum_{\x} | \sum_{\h \in \Hcal} w_B(\x+
q_2\h)\chi_{q_1}(f(\x+ q_2\h ))|^2 =
\sum_{\h_1 \in \Hcal}\sum_{\h_2 \in \Hcal}S(\h_1,\h_2)\]
where
\[S(\h_1,\h_2)=\sum_{\x}\chi_{q_1}(f(\x+ q_2\h_1))
\overline{\chi_{q_1}(f(\x+ q_2\h_2 ))}w_B(\x+ q_2\h_1) w_B(\x+ q_2\h_2 ).\]
We then see that $S(\h_1,\h_2)=S(\h_1-\h_2,\mathbf{0})$, and hence
that
\begin{eqnarray*}
\Sig_2 &\leq & \sum_{\h_1\in\Hcal}\sum_{\h_2\in\Hcal}
S(\h_1-\h_2,\mathbf{0}) \\
& =&  \sum_{\h \in \Hcal_0} \prod_{j=1}^{n}(H-|h_j|)
S(\h,\mathbf{0})\\
&\le &H^n \sum_{\h \in \Hcal_0} | \sum_{\x} \chi_{q_1}(f(\x+ q_2\h))
\overline{\chi_{q_1}(f(\x))}  w_B(\x+ q_2\h) w_B(\x)|,
\end{eqnarray*}
where $\Hcal_0 = [-H,H]^n$.
We now further split $\Sig_2$ into the single term with $\h =
(0,\ldots,0)$, which we will call $\Sig_{2A}$, and the remainder of
the sum over $\h \in \Hcal_0$, $\h \neq 0$, which we will call
$\Sig_{2B}$.

Our goal is to give upper bounds for $\Sig_1$ and $\Sig_2$. The first
admits a trivial bound: clearly, \beq\label{Sig1_bd}
\Sig_1 \ll B^n.
\eeq
Note that as we apply a trivial bound to this term, which is the only sum whose modulus is $q_2$, we do not need to assume that $q_2$ is square-free. (This is what enables us to treat $S(\U)$
as part of the main sieve in Lemma \ref{sieve_lemma}, but not $S(\V)$.)

Similarly, we bound $\Sig_{2A}$ trivially as
\beq\label{Sig2A_bd}
\Sig_{2A} \leq H^n\sum_{\x \in \Z^n} (w_B(x))^2 \ll H^nB^n.
\eeq

Combining (\ref{Sig1_bd}) and (\ref{Sig2A_bd}) in (\ref{H2T2}), we
have now shown that \[T(q_1,q_2)\ll H^{-n}\Sig_1^{1/2}(\Sig_{2A} + \Sig_{2B})^{1/2}
\ll H^{-n}B^{n/2}(H^nB^n + \Sig_{2B})^{1/2},\]
whence
\beq\label{Tq_upper}
T(q_1,q_2)\ll B^{n/2}q_2^{n/2}+B^{-n/2}q_2^n\Sig_{2B}^{1/2}.
\eeq

We now require a nontrivial upper bound
for $\Sig_{2B}$.  We may write
\beq\label{Sig_T}
\Sig_{2B} \ll H^n \sum_{\bstack{\h \in \Hcal_0}{\h \neq 0}} |T(\hbf)|,
	\eeq
where \beq\label{T_triple}
T(\hbf) = \sum_{k \modd{q_1}}\chi_{q_1}(k)
\sum_{\tstack{\xbf \in \Z^n}{f(\xbf + q_2\hbf)-kf(\x)\con 0 \modd{q_1}}
{(f(\xbf),q_1)=1}} w_{B,\hbf}(x).
\eeq
Here we have set $w_{B,\hbf}(x) = w_B(\x + q_2\hbf)w_B(\x)$.

In order to bound $T(\hbf)$, we shall consider a general sum of the form
\[S=\sum_{\bstack{\x\in\Z^n}{q\mid h(\x),\,(g(\x),q)=1}}W(\x/L)\]
in which $q$ is either prime or a product $p_1p_2$ of primes satisfying $p_1<p_2<2p_1$. We suppose that $h(\x)$ and $g(\x)$ are
integral polynomials in $\x=(x_1,\ldots,x_n)$, with ${\rm deg}(h)\le
{\rm deg}(g)=d$, where $d\ge 3$.  We shall take the
weight function $W(\x)$ to be smooth and supported on $[-1,1]^n$, and we shall write $\Delta$ for
the maximum of the moduli of all partial derivatives
of $W$ with order at most $n+1$.

Under these assumptions we shall estimate $S$, using information on the behaviour of $h(\x)$ and $g(\x)$ modulo the prime factors 
of $q$.  Let $G(\x)$ be the leading form for $g(\x)$, so that $G(\x)$
has degree $d$.  We will require $G(\x)$ to be nonsingular modulo every
prime factor $p$ of $q$.  We shall assume further that either the leading form for $h(\x)$ is a constant multiple of $G(\x)$ or that the degree of $h(\x)$ is strictly less than $d$.  It follows that there is exactly one value $\gamma$ modulo $p$ for which $h(\x)-\gamma g(\x)$ has
degree less than $d$, when considered over $\F_p$.  If
$H(\x)\in\F_p[\x]$ is the leading form for $h(\x)-\gamma g(\x)$ we shall
require $H$ to have degree at least 2, and we write $s(h,g;p)$ for the
dimension of the singular locus of the variety $H(\x)=0$ in $\mathbb{A}^n(\mathbb{F}_p)$.

In the situation above, the leading form of $ah(\x)+bg(\x)$ will
have a singular locus of dimension at most $s(h,g;p)$, for any
$(a,b)\in\F_p^2-\{(0,0)\}$.  In particular, if $s(h,g;p)=0$ it follows
from the fundamental theorem of Deligne \cite{Del74} that
\beq\label{deligne}
\sum_{\x\modd{p}}e_p(ah(\x)+bg(\x)+\mathbf{v}\cdot \x)\ll_{n,d}p^{n/2}.
\eeq

Using this bound we shall ultimately establish in Section \ref{sec_sums_div} the following result.
\begin{prop}\label{Lemma_Poisson_fg}
Adopt the assumptions above, and let $p=q$ if $q$ is prime, or $p=p_1$
if $q=p_1p_2$.  Then if $L\ge 1$ and $n\ge 2$ we have
\begin{eqnarray}\label{lemma_sum}
S&=&q^{-2}\phi(q)\sum_{\x\in\Z^n}W(\x/L)
+O_{n,d}(\Delta L^sq^{(n-s)/2})\nonumber\\
&&\hspace{2cm}\mbox{}+O_{n,d}(\Delta L^n p^{(s-n+2)/2}q^{-1}).
\end{eqnarray}
Here we have set $s=s(h,g;q)$ if $q$ is
prime, or \[s=\min\left(s(h,g;p_1)\,,\,s(h,g;p_2)\right)\]
if $q=p_1p_2$.
\end{prop}

In our application $\Delta$ will be $O_{n,d}(1)$.  However since our
proof of Proposition \ref{Lemma_Poisson_fg} uses an induction in which
the weight $W$ varies, we have found it clearer to include $\Delta$ in
the error estimates above.

We apply Proposition \ref{Lemma_Poisson_fg} to the innermost sum in
(\ref{T_triple}) with $q=q_1$, \[h(\xbf) = f(\xbf+q_2\hbf)-kf(\x),\qquad g(\xbf) = f(\xbf),\]
$W(\xbf)=w(\xbf+q_2\hbf)w(\xbf)$,
and $L=B$. Note that $W(\xbf)$ is then supported on a cube of side 2,
and so $\Delta\ll_{n,d} 1$. Note also that the condition that the leading form $G(\xbf)$ of $g(\xbf)$ is nonsingular modulo every prime factor of $q_1$ is satisfied, provided that $B \gg1$. Conveniently, since the main term in (\ref{lemma_sum}) is independent
of $k$, its total contribution to (\ref{T_triple}) when summed over $k$ is zero. In order to
estimate $\Sigma_{2B}$ via (\ref{Sig_T}) we need to understand how
 \[s=s(h,g;p) = s(f(\xbf+q_2\hbf)-kf(\x),f(\x);p)\]
varies as we change $\hbf$.
The leading form of $f(\xbf+q_2\hbf)-kf(\x)-\gamma f(\x)$, taken over
$\F_p$, can only have degree less than $d$ in the case $k+\gamma=1$, in
which case the terms of degree $d-1$ are $q_2\hbf\cdot\nabla
F(\x)$, where $F(\x)$ is the leading form of $f$. Thus we may interpret $s(h,g;p)$ as the dimension of the
singular locus of the variety $\hbf\cdot\nabla F(\x)=0$ in $\mathbb{A}^n(\F_p)$.
Our next lemma provides the necessary information about this.

\begin{lemma}\label{s,H est}
Suppose that $F(\x)\in\F_p[x_1,\ldots,x_n]$ is a nonsingular form
of degree $d$, and let $H$ be a positive integer.
Then if $0\le s\le n$, the number of non-zero $\h\in [-H,H]^n$ for which
the variety $\hbf\cdot\nabla F(\x)=0$ has singular locus of affine
dimension $s$ is $O_{n,d}(H^{n-s}+H^np^{-s})$.
\end{lemma}
We will prove this in \S \ref{last}.
It follows immediately from Lemma \ref{s,H est} that the number of non-zero $\hbf\in\Hcal_0$ for which
$s(h,g;q_1)=s$ will be $O_{n,d}(H^{n-s}+H^nq_1^{-s})$ if $q_1$ is a prime.  On the other hand, if $q_1=p_1p_2$ with $p_1<p_2 <2 p_1$ then the
number of $\hbf$ with $\min(s(h,g;p_1),s(h,g;p_2))=s$ will be $O_{n,d}(H^{n-s}+H^np^{-s})$ where $p=p_1$ or $p_2$. Thus in either case we may write the bound as
$O_{n,d}(H^{n-s}+H^np^{-s})$, in the notation of Proposition \ref{Lemma_Poisson_fg}. The error terms in (\ref{lemma_sum}) therefore contribute to (\ref{Sig_T}) a total of
\[\ll_{n,d}q_1H^n\sum_{0\le s\le n}(H^{n-s}+H^np^{-s})
(B^sq_1^{(n-s)/2}+B^np^{(s-n+2)/2}q_1^{-1}).\]
Each summand takes the form $XY^s$ as a function of $s$ and is
therefore maximal either at $s=0$ or $s=n$.  From $s=0$ we get a
contribution
\begin{eqnarray*}
&\ll_{n,d}&q_1H^n(H^{n}+H^n)(q_1^{n/2}+B^np^{-(n-2)/2}q_1^{-1})\\
&\ll_{n,d}&H^{2n}q_1^{(n+2)/2}+H^{2n}B^np^{-(n-2)/2},
\end{eqnarray*}
while for $s=n$ we obtain
\begin{eqnarray*}
&\ll_{n,d}&q_1H^n(1+H^np^{-n})(B^n+B^npq_1^{-1})\\
&\ll_{n,d}& H^nB^nq_1+H^{2n}B^np^{-n}q_1.
\end{eqnarray*}
We therefore conclude that an overall bound for $\Sig_{2B}$ in (\ref{Sig_T}) is
\[\Sig_{2B} \ll_{n,d}H^{2n}q_1^{(n+2)/2}+H^{2n}B^np^{-(n-2)/2}+H^nB^nq_1.\] Proposition \ref{T_prop} now follows from (\ref{Tq_upper}) on
recalling that $H=[B/q_2]$.

\section{Bounding the sieve terms}\label{bst}
We are now ready to apply Proposition \ref{T_prop} to bound the main
sieve and the prime sieve in (\ref{N_upper_bd}) and prove Theorem \ref{form_thm}.
 The main sieve is
bounded above by:
\[A^{-2}\sum_{v,v'\in\V}\sum_{u\neq u'\in\U}
\left|\sum_n\om(n)\chi_{uv}(n)\overline{\chi_{u'v'}(n)}\right| \ll\sup_{v,v'\in\Vcal}\sup_{u\neq u'\in \Ucal} |T(uu',vv')| ,\]
where $T(uu',vv')$ is defined as in (\ref{Tq_dfn}) with $q_1=uu'$, $q_2=vv'$, and the characters
$\chi_{q_1}^*, \chi_{q_2}^*$ as defined in (\ref{chi_dfn}).  According
to the definitions (\ref{Udef}) and (\ref{Vdef}) for the sieving sets $U$ and $V$, Proposition \ref{T_prop} shows that the
above is
\[\ll_{f} B^{n/2}Q^{n-(n-1)\alpha}+B^{n/2}Q^{(n+2)\alpha/2}
+B^nQ^{-(n-2)\alpha/4}.\]
We choose $\alpha=2/3$ so as to match the first two terms above,
giving a bound
\beq\label{sieve1}
\ll_f B^{n/2}Q^{(n+2)/3}+B^nQ^{-(n-2)/6}
\eeq
for the main sieve term.  This is subject to the condition $q_2\le B$,
for which it suffices to have $4Q^{2/3}\le B$.

We now turn to the prime sieve, given in (\ref{sig_lemma_eqn}) as
\begin{eqnarray*}
UA^{-2} \sum_{v \neq v' \in \V} \left| \sum_{n} \om(n) \chi_v(n)
 \overline{\chi_{v'}(n)} \right|
	&\ll& UV^2A^{-2} \sup_{v \neq v' \in \Vcal} |T(v,v')|\\
&\ll& U^{-1} \sup_{v \neq v' \in \Vcal} |T(v,v')|,
\end{eqnarray*}
where $T(v,v')$ is again defined as in (\ref{Tq_dfn}) but with respect
to characters $\chi_v, \overline{\chi_{v'}}$ with prime moduli. Since $v$ and
$v'$ are each of order $Q^{1-\alpha}=Q^{1/3}$ we get the immediate bound
\begin{eqnarray*}
&\ll_{n,d}&U^{-1}\{B^{n/2}Q^{(n+1)/6}+B^{n/2}Q^{(n+2)/12}+B^nQ^{-(n-2)/12}\}\\
&\ll_{n,d}&Q^{-2/3}(\log Q)\{B^{n/2}Q^{(n+1)/6}+B^nQ^{-(n-2)/12}\}.
\end{eqnarray*}
On combining this bound with (\ref{lead_term}) and (\ref{sieve1}) and inserting the result
into Lemma \ref{sieve_lemma}, we find that
\begin{eqnarray*}
\sum_{n\neq 0}\om (n^{r})
	&\ll_f & B^nQ^{-1}(\log Q)^2 +B^{n/2}Q^{(n+2)/3}+B^nQ^{-(n-2)/6} \\[-12pt]
		&& \qquad + \; Q^{-2/3}(\log Q)\{B^{n/2}Q^{(n+1)/6}+B^nQ^{-(n-2)/12}\}\\
&\ll_f & (\log Q)^2\{B^n(Q^{-1}+Q^{-(n-2)/6} + Q^{-2/3 -(n-2)/12})+B^{n/2}Q^{(n+2)/3}\}.
\end{eqnarray*}
The optimal choice of $Q$ will be
\[Q=\left\{\begin{array}{cc} B^{3n/(2n+10)}, & n\ge 8,\\
B^{3n/(3n+2)}, & 2\le n\le 8,\end{array}\right.\]
yielding bounds
\[\ll_f B^{n-3n/(2n+10)}(\log B)^2\]
and
\[\ll_f B^{n-n(n-2)/(6n+4)}(\log B)^2,\]
respectively. Theorem \ref{form_thm} then follows.

\section{Proof of Proposition \ref{Lemma_Poisson_fg}}
\label{sec_sums_div}
Our treatment of Proposition \ref{Lemma_Poisson_fg}, which is essentially a version of Poisson summation, is motivated by the
argument used by Heath-Brown \cite[Theorem 3]{HB94}, and employs
induction on $s$.
We therefore begin by establishing the base case for the induction, in
which $s=0$. We split the values of $\x$ into residue classes modulo
$q$ and use the Poisson Summation Formula to obtain
\begin{eqnarray}\label{PS}
S&=&\sum_{\bstack{\zbf\modd{q}}{q|h(\zbf),\,(g(\zbf),q)=1}}
\sum_{\ubf \in \Z^n} W\left(\frac{\zbf + q\ubf}{L}\right) \nonumber\\
& = & \left(\frac{L}{q}\right)^n \sum_{\vbf \in \Z^n}
       \widehat{W}\left(\frac{L\vbf}{q}\right) S_q(\vbf),
	\end{eqnarray}
where
\[ S_q(\vbf) = \sum_{\bstack{\zbf\modd{q}}{q|h(\zbf),\,(g(\zbf),q)=1}}
e_q(\vbf \cdot \zbf).\]
We may estimate
\[\widehat{W}(\x)=\int_{\R^n}W(\y)e(-\x\cdot \y)d\y\]
by integrating by parts $n+1$ times with respect to $y_j$, say. This shows that $\widehat{W}(\x)\ll_n\Delta |x_j|^{-n-1}$, and since $j$ is
arbitrary we may conclude that
\beq\label{whest}
\widehat{W}(\x)\ll_n\Delta|\x|^{-n-1}, \qquad \text{for $|x| \geq 1$};
\eeq
for $|x| \leq 1$, we will employ the trivial bound $\widehat{W}(\x) \ll_n \Delta$.

When $q=p_1p_2$ the sum $S_q(\vbf)$ satisfies a multiplicativity relation
\[S_{q}(\vbf)=S_{p_1}(\vbf)S_{p_2}(\vbf).\]
Moreover if $p$ is prime then
\begin{eqnarray*}
S_p(\vbf)&=&\sum_{p|h(\zbf)}e_q(\vbf \cdot \zbf)-
\sum_{p|h(\zbf),g(\zbf)}e_q(\vbf \cdot \zbf)\\
&=&p^{-1}\sum_{a\modd{p}}\sum_{\zbf\modd{p}}e_q(ah(\zbf)+\vbf\cdot\zbf)\\
&&\hspace{1cm}\mbox{}-p^{-2}\sum_{a,b\modd{p}}\sum_{\zbf\modd{p}}
e_q(ah(\zbf)+bg(\zbf)+\vbf\cdot\zbf).
	\end{eqnarray*}
Since we are
assuming that $s=0$, the Deligne estimate (\ref{deligne}) applies when $a\not=0$, for the first
sum above, and for $(a,b)\not=(0,0)$ for the second.  We therefore have
\[S_p(\vbf)=(p^{-1}-p^{-2})\sum_{\zbf\modd{p}}e_p(\vbf\cdot\zbf)
+O_{n,d}(p^{n/2}).\]
It follows that \beq\label{s0}
S_p(\mathbf{0})=\phi(p)p^{n-2}+O_{n,d}(p^{n/2})
\eeq
and that
$S_p(\vbf)=O_{n,d}(p^{n/2})$ for $p\nmid\vbf$. Using
the multiplicativity relation, we now see by (\ref{whest})
that terms in (\ref{PS}) with $\vbf \neq 0$, $\vbf$
coprime to $q$ contribute
\begin{eqnarray*}
&\ll_{n,d}&\left(\frac{L}{q}\right)^n q^{n/2}\sum_{\vbf\in\Z^n-\{\mathbf{0}\}}
\left|\widehat{W}\left(\frac{L\vbf}{q}\right)\right|\\
&\ll_{n,d}& L^n q^{-n/2}\sum_{\vbf \in \Z^n-\{\mathbf{0}\}}
\Delta \min \{ 1, \left(\frac{q}{L|\vbf|}\right)^{n+1}\}\\
&\ll_{n,d}& L^n q^{-n/2}\Delta\left(\frac{q}{L}\right)^n\\
&\ll_{n,d}& \Delta q^{n/2}.
\end{eqnarray*}
If $q=p_1p_2$ then the terms with $p_1\mid\vbf$ but $p_2\nmid\vbf$
have $S_q(\vbf)\ll q^n$ and hence contribute
\begin{eqnarray*}
&\ll_{n,d}&\left(\frac{L}{q}\right)^n q^n
\sum_{\bstack{\vbf \in \Z^n-\{\mathbf{0}\}}{p_1|\vbf}}
\left|\widehat{W}\left(\frac{L\vbf}{q}\right)\right|\\
&\ll_{n,d}& L^n\sum_{\ubf \in \Z^n-\{\mathbf{0}\}}
\Delta \min \{ 1, \left(\frac{q}{Lp_1|\ubf|}\right)^{n+1} \}\\
&\ll_{n,d}& L^n\Delta \left(\frac{q}{p_1 L}\right)^n \\
&\ll_{n,d}& \Delta q^{n/2},
\end{eqnarray*}
and similarly if $p_2\mid\vbf$ but $p_1\nmid\vbf$. We therefore deduce that
\[S=\left(\frac{L}{q}\right)^n S_q(\mathbf{0})
\sum_{\bstack{\vbf \in \Z^n}{q|\vbf}}
 \widehat{W}\left(\frac{L\vbf}{q}\right) +O_{n,d}(\Delta q^{n/2}).\]
However (\ref{s0}) yields
$S_q(\mathbf{0})=\phi(q^{n-1})+O_{n,d}(q^{n/2})$ if $q$ is
prime and similarly
$S_q(\mathbf{0})=\phi(q^{n-1})+O_{n,d}(q^{(3n-2)/4})$ if $q=p_1p_2$.  Moreover
\[\sum_{\bstack{\vbf \in \Z^n}{q|\vbf}}\widehat{W}\left(\frac{L\vbf}{q}\right) =\sum_{\ubf \in \Z^n}\widehat{W}(L\ubf) =L^{-n}\sum_{\ubf \in \Z^n}W\left(L^{-1}\ubf\right)\ll\Delta.\]
The case $s=0$ of the proposition then follows.

When $n=2$ and $s=1$ or $2$, the proposition is immediate. To see this we
observe that the polynomial $h(\xbf)$ cannot vanish identically modulo a prime divisor $p$ of $q$, by
our initial assumption that $h(\x)-\gamma g(\x)$ has degree at least 2,
but strictly less than $d$.  We then estimate $S$ via the following
lemma.
\begin{lemma}\label{bhb}
Suppose that $q$ is either prime or the product $p_1p_2$ of primes
$p_1<p_2<2p_1$. Suppose $k\le n$ and that for $p=q$ (in the first case) or for
$p=p_1$ and $p=p_2$ (in the second) we are given a variety
$V_p\subseteq\mathbb{A}^n(\mathbb{F}_p)$ of dimension $k$ and degree at
most $D$.  Then if $\tau_p$ is the natural map from $\Z$ to $\F_p$,
we have
\[\#\{\x\in\Z^n\cap[-R,R]^n:\tau_p(\x)\in V_p\mbox{ for }p\mid q\}
\ll_{n,D} R^n q^{k-n}+R^k.\]
\end{lemma}
This is a special case of Lemma 4 of Browning and
Heath-Brown's work \cite{BHB09}, in which we take
$W=\mathbb{A}^n(\Q)$, $l=n$, and $k_i=k$ in their notation.

We apply the lemma with $n=2$ and $k=1$ to give
$S\ll_{n,d}L^2q^{-1}+L$ in our situation. We also have
\[q^{-n}\phi(q^{n-1})\sum_{\x\in\Z^n}W(\x/L)\ll L^2q^{-1}.\]
We therefore see that these are dominated by the error terms in
(\ref{lemma_sum}) if $n=2$ and $s=1$ or $2$.

We turn now to the induction argument, for which we assume $n\ge 3$
and $s\ge 1$. The induction step
will reduce both $n$ and $s$ by 1, giving us a case for which we
already know that the proposition holds.  The plan is to choose a suitable
matrix $M\in{\rm SL}_n(\Z)$, and to work with polynomials $h_{M,c}(\y)$
and $g_{M,c}(\y)$ in $n-1$ variables $\y=(y_1,\ldots,y_{n-1})$ defined
by setting
\[h_M(\x)=h(M\x)\qquad\mbox{and}\qquad h_{M,c}(\y)=h_M(\y,c),\]
and similarly for $g$. If we also set
\[W_M(\x)=W(M\x)\qquad\mbox{and}\qquad W_{M,c}(\y)=W_M(\y,c)\]
we then find that
\[S=\sum_{c\in\Z}\sum_{\bstack{\y\in\Z^{n-1}}
{q|h_{M,c}(\y),\,(g_{M,c}(\xbf),q)=1}}W_{M,c}(\y/L).\]
In order to apply the induction hypothesis we use the following lemma
to provide a suitable matrix $M$.

\begin{lemma}\label{M}
Suppose that $q$ and the polynomials $h$ and $g$ are as in the preamble
to Proposition \ref{Lemma_Poisson_fg}. Then there is a matrix $M\in{\rm
 SL}_n(\Z)$ with entries bounded in modulus by $||M||\ll_{n,d}1$, and
having the following properties for every prime divisor $p$ of $q$.
Firstly, the leading form for $g_{M,c}$ will be nonsingular modulo $p$,
and secondly, the leading form for $h_{M,c}-\gamma g_{M,c}$ will have degree at least 2 over $\F_p$, with
singular locus of dimension at most $\max(s(h,g;p)-1,0)$.
\end{lemma}

We will prove this in the next section, but we first show how we can then complete
the induction step. We first note that
$||M^{-1}||\ll_n||M||^{n-1}\ll_{n,d} 1$.  Hence if $W_{M}(\x)\not=0$
we have $M\x\ll_n 1$, and hence $\x\ll_{n,d}1$.
It follows that $W_{M,c}(\y/L)$ vanishes unless $c\ll_{n,d}L$ and that $W_{M,c}(\tbf)$ has support $\tbf\in [-c_0,c_0]^{n-1}$ with $c_0\ll_{n,d}1$. We therefore write
$W_0(\tbf)=W_{M,c}(c_0\tbf)$ so that $W_0(\tbf)$ is supported in
$[-1,1]^{n-1}$.  We also observe that any $j$-th order partial
derivative of $W_0$ is of size $O_{n,d}(\Delta)$. We may now apply Lemma \ref{Lemma_Poisson_fg} with
$s$ replaced by $s-1$ to find that
\begin{eqnarray*}
\lefteqn{\sum_{\bstack{\y\in\Z^{n-1}}
{q|h_{M,c}(\y),\,(g_{M,c}(\xbf),q)=1}}W_{M,c}(\y/L)}\\
&=&\sum_{\bstack{\y\in\Z^{n-1}}
{q|h_{M,c}(\y),\,(g_{M,c}(\xbf),q)=1}}W_0(c_0^{-1}L^{-1}\y)\\
&=&q^{1-n}\phi(q^{n-2})\sum_{\y\in\Z^{n-1}}W_0(c_0^{-1}L^{-1}\y)
+O_{n,d}(\Delta L^{s-1}q^{((n-1)-(s-1))/2})\\
&&\hspace{2cm}\mbox{}+O_{n,d}(\Delta L^{n-1}p^{((s-1)-(n-1)+2)/2}q^{-1}).
\end{eqnarray*}
When we sum over all $c$ such that $W_{M,c}(\ybf/L) \neq 0$,  the error terms contribute
\begin{eqnarray*}
&\ll_{n,d}&L\Delta L^{s-1}q^{((n-1)-(s-1))/2}
+L\Delta L^{n-1}p^{((s-1)-(n-1)+2)/2}q^{-1}\\
&\ll_{n,d}&\Delta L^sq^{(n-s)/2}+\Delta L^np^{(s-n+2)/2}q^{-1}.
\end{eqnarray*}
Moreover we have
\begin{eqnarray*}
\sum_{c\in \Z} \sum_{\y\in\Z^{n-1}}W_0(c_0^{-1}L^{-1}\y)&=&
\sum_{c \in \Z} \sum_{\y\in\Z^{n-1}}W_{M,c}(\y/L)\\
&=&\sum_{\x\in\Z^{n}}W(\x/L).
\end{eqnarray*}
It will then follow that
\[S=q^{1-n}\phi(q^{n-2})\sum_{\x\in\Z^{n}}W(\x/L)
+O_{n,d}(\Delta L^sq^{(n-s)/2})+O_{n,d}(\Delta L^np^{(s-n+2)/2}q^{-1}),\]
which suffices for our induction step.

\section{Proof of Lemmas \ref{s,H est} and \ref{M}}\label{last}

Our proof of Lemma \ref{s,H est} is based on the following result of
Heath-Brown \cite[Lemma 2]{HB94}
\begin{lemma}\label{nab}
Let $F(\xbf)\in\F_p[x_1,\ldots,x_n]$ be a smooth form of degree $d$. For each $\hbf \in \overline{\F_p}^n$, let $S_\hbf$ denote the affine variety
\[ S_\hbf = \{ \x : \hbf \cdot \nabla^2 F(\xbf)=0 \},\]
and for every non-negative integer $s \leq n$, let
\[ T_s= \{ \hbf: \dim(S_\hbf(F)) \geq s \}.\]
Then $T_s$ is an affine variety, and has dimension at most $n-s$.  Moreover it may be defined by $O_{n,d}(1)$ equations, each of degree $O_{n,d}(1)$.
\end{lemma}

Clearly Lemma \ref{s,H est} follows from this estimate in conjunction
with Lemma \ref{bhb}.

We turn now to the proof of Lemma \ref{M}.  We recall that $G$
is the leading form of $g$ (and is assumed to be nonsingular modulo every prime divisor of $q$) and that $H$ is the leading form of
$h-\gamma g$.  Thus the leading form of $g_{M,c}$ will be $G_{M,0}$,
and similarly the leading form of $h_{M,c}-\gamma g_{M,c}$ will be $H_{M,0}$, providing that $G_{M,0}$ and $H_{M,0}$ do not
vanish identically.  We may view the variety in $\mathbb{A}^n(\F_p)$
defined by $G_{M,0}(\ybf)=0$ as being the intersection of $G_M(\x)=0$
with the hyperplane $x_n=0$.  This is isomorphic to the intersection
of the variety
\[\mathcal{G}_p:\;G(\x)=0\]
with $(M^{-1}\x)_n=0$.  Thus if $\mbf$ is the column
vector whose transpose is the bottom row of $M^{-1}$, the variety in
which we are interested will be \[\mathcal{G}_p^{\mbf}:\; G(\x)=\mbf\cdot\x=0.\]
It will be convenient to use the notation $s(V)$ for the affine dimension of the singular locus of a variety $V$.  Thus to confirm the first conclusion of Lemma \ref{M}, we are
hoping to show that $s(\mathcal{G}_p^{\mbf})=0$, and hence $g_{M,c}$ is nonsingular modulo $p$, for a suitable matrix
$M$.

We now recall Lemma 5 of Heath-Brown \cite{HB94}, which states
that for any prime $p$ and any form $R(\x)\in\F_p[x_1,\ldots,x_n]$ one
has $s(\mathcal{R}_p^{\mbf})\geq s(\mathcal{R}_p)-1$ for all non-zero
$\mbf\in\F_p^n$, where $\mathcal{R}_p$ and $\mathcal{R}_p^{\mbf}$ are
defined analogously to the case for $G$ above.  Moreover there
exists a non-zero form $\widehat{R}_p$ depending on $p$ and $R$ such that the
degree of $\widehat{R}_p$ is bounded in terms of $n$ and the degree of $R$ alone, and such that \[s(\mathcal{R}_p^{\mbf})=\max\left(s(\mathcal{R}_p)-1\,,\,0\right)\] whenever $p \ndiv \widehat{R}_p(\mbf)$.

Thus in our case, if $p\ndiv\widehat{G}_p(\mbf)$ then $\mathcal{G}_p^{\mbf}$ will be
nonsingular, since $s(G)=0$ and so $s(\mathcal{G}_p) \leq 1$. In exactly the same way we find that if $p\ndiv\widehat{H}_p(\mbf)$ then \[s(\mathcal{H}_p^{\mbf})=\max\left(s(\mathcal{H}_p)-1\,,\,0\right),\] and in particular $H_{M,0}$ will not vanish identically.

We therefore wish to find a vector $\mbf$ such that
$q\ndiv\widehat{G}_q(\mbf)\widehat{H}_q(\mbf)$, if $q$ is prime, or such that
$p\ndiv\widehat{G}_p(\mbf)\widehat{H}_p(\mbf)$ for $p=p_1$ and $p=p_2$ in the
case $q=p_1p_2$.  However, according to Lemma \ref{bhb},
if one has a nonzero polynomial
$f(\x)\in\F_p[x_1,\ldots,x_n]$ of degree $D$, then
\[\#\{\x\in(0,T]^n:\,f(\x)\con 0\modd{p}\}\ll_{n,D} T^n p^{-1} + T^{n-1}.\]
In our case we deduce that, if $T\gg_{n,d}1$ and $q,p_1,p_2\gg_{n,d}1$, then there
will be a vector $\mbf\in(0,T]^n$, such that none of
$q\mid\widehat{G}_q(\mbf)\widehat{H}_q(\mbf)$ or
$p\mid\widehat{G}_p(\mbf)\widehat{H}_p(\mbf)$ holds.  Clearly we may suppose
that $\mbf$ is primitive, since we can divide out by any common factor
without affecting the non-divisibility result.  Proposition
\ref{Lemma_Poisson_fg} is of course trivial if $q\ll_{n,d}1$ and so we
may therefore conclude that there is an admissible primitive $\mbf\ll_{n,d}1$.

Finally, to finish the proof of Lemma \ref{M} we
observe that given such a vector $\mbf$ there is a matrix $M_1\in{\rm
 SL}_n(\mathbb{Z})$ having the transpose of $\mbf$ as its last row,
and such that $||M_1||\ll_{n,d}$.  We then find that $M$ defined by $M^{-1}= M_1$
is acceptable for Lemma \ref{M}. This completes the proof of Lemma
\ref{M}.

\section{Acknowledgements}

Part of this work was carried out while Roger Heath-Brown was visiting
the Mathematical Sciences Research Institute, Berkeley.  The
hospitality and financial support of the institute is gratefully
acknowledged. Lillian Pierce was supported by a Marie Curie Fellowship funded by the European Commission for the duration of this work.

\bibliographystyle{amsplain}
\bibliography{NoThBibliography}

\end{document}